\documentclass[letterpaper, 10 pt, conference]{ieeeconf}  

\IEEEoverridecommandlockouts                              

\usepackage{graphics, subfig} 
\usepackage{epsfig} 
\usepackage{mathptmx} 
\usepackage{times} 
\usepackage{amsmath} 
\usepackage{amssymb}  
\usepackage{algorithmic,algorithm}
\usepackage{comment}
\usepackage{xcolor}
\newtheorem{theorem}{Theorem}
\newtheorem{remark}{Remark}

{\renewcommand{\baselinestretch}{1}}

\begin{document}

\title{\LARGE \bf
Multiway $k$-Cut in Static and Dynamic Graphs: A Maximum \\Entropy Principle Approach
}
\author{Mayank Baranwal$^{1}$, Amber Srivastava$^{2a}$, and Srinivasa M. Salapaka$^{2b}$
\thanks{$^{1}$Department of Electrical and Computer Engineering at the University of Michigan, Ann-Arbor, 48109 MI, USA. {\tt\small mayankb@umich.edu}}
\thanks{$^{2}$Department of Mechanical Science and Engineering, University of Illinois at Urbana-Champaign, 61801 IL, USA. $^{a}${\tt\small asrvstv6@illinois.edu}, $^{b}${\tt\small salapaka@illinois.edu}}
\thanks{The authors would like to acknowledge NSF grants ECCS 15-09302, CMMI 14-63239 and CNS 15-44635 for supporting this work.}}

\thispagestyle{empty}
\pagestyle{empty}

\maketitle

\begin{abstract}
This work presents a maximum entropy principle based algorithm for solving minimum multiway $k$-cut problem defined over static and dynamic {\em digraphs}. A multiway $k$-cut problem requires partitioning the set of nodes in a graph into $k$ subsets, such that each subset contains one prespecified node, and the corresponding total cut weight is minimized. These problems arise in many applications and are computationally complex (NP-hard). In the static setting this article presents an approach that uses a relaxed multiway $k$-cut cost function; we show that the resulting algorithm converges to a local minimum. This iterative algorithm is designed to avoid poor local minima with its run-time complexity as $\sim O(kIN^3)$, where $N$ is the number of vertices and $I$ is the number of iterations. In the dynamic setting, the edge-weight matrix has an associated dynamics with some of the edges in the graph capable of being influenced by an external input. The objective is to design the dynamics of the controllable edges so that multiway $k$-cut value remains small (or decreases) as the graph evolves under the dynamics. Also it is required  to determine the time-varying partition that defines the minimum multiway $k$-cut value. Our approach is to choose a relaxation of multiway $k$-cut value, derived using maximum entropy principle, and treat it as a control Lyapunov  function to design control laws that affect the weight dynamics. Simulations on practical examples of interactive foreground-background segmentation, minimum multiway $k$-cut optimization for non-planar graphs and dynamically evolving graphs that demonstrate the efficacy of the algorithm, are presented.
\end{abstract}

\section{Introduction}\label{sec:Intro}

\begin{figure}
	\centering
	\begin{tabular}{c}
		\includegraphics[width=\columnwidth]{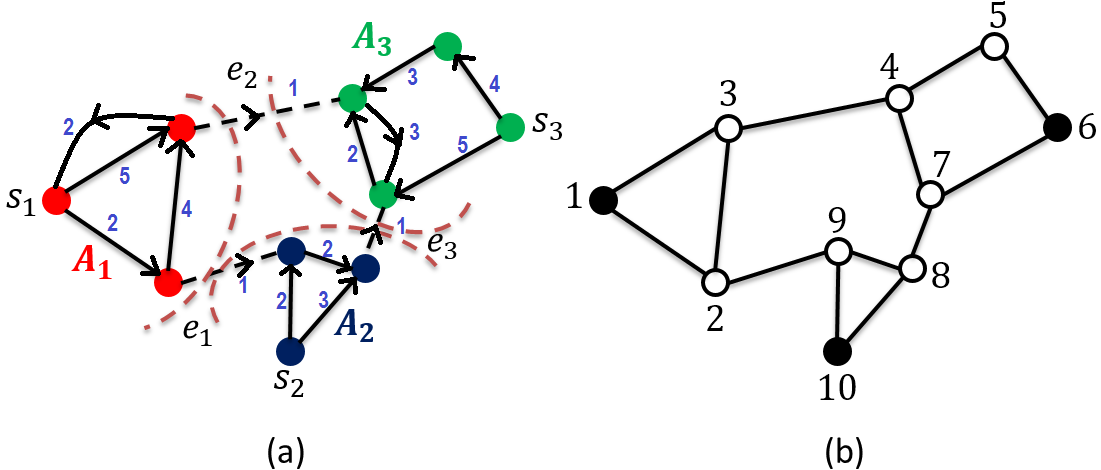}
	\end{tabular}
	\caption{{\small (a) Schematic of a minimum multiway $k$-cut problem with $k=3$ and $S=\{s_1,s_2,s_3\}$. Here $A_1,A_2$ and $A_3$ are the disconnected components. (b) Illustrative example of a multiway $3$-cut problem with $S=\{1,6,10\}$ and unit edge-weights described in Sec. \ref{sec:sim_results}.}}
	\label{fig:multiwayCut}
\end{figure}

The {\em Multiway $k$-Cut} Problem \cite{lowen2011multiway} is a generalization of minimum {\em s-t} cut problem and has applications in parallel and distributed computing \cite{stone1977multiprocessor}, as well as in chip design. Multiway cut also finds applications in several other problems of related interest such as extending a partial $k$-coloring of a graph \cite{erdHos1994weighted}. Given a graph $G=(V,E)$ with vertex set $V, |V| = N\in\mathbb{N}$, set $E$ of edges, edge weights $w:E\to\mathbb{R}^+$, and a set of terminals $S = \{s_1,s_2,\dots,s_k\}\subseteq V$, a {\em multiway $k$-cut} is a set of edges whose removal disconnects the terminals from each other. The goal of the minimum multiway $k$-cut problem is to find a minimum weight set $E'\subseteq E$ of edges such that removing $E'$ from $G$ separates all terminals. Fig. \ref{fig:multiwayCut}a shows a schematic of a minimum multiway $3$-cut problem with set of terminals denoted by $S=\{s_1,s_2,s_3\}$. The objective is to obtain a partition of the vertex set $V$ into disconnected components $\{A_j:A_j\subset V,1\leq j\leq 3\}$, such that $s_j\in A_j$ for all $j$, and the total cut size, $\frac{1}{2}\sum_{j=1}^{3}w(A_j,\bar{A_j})$ is minimized.

While the problem of computing a minimum {\em s-t} cut (i.e., $k=2$) in static graphs is solvable in polynomial time, it is shown in \cite{dahlhaus1994complexity} that the minimum multiway $k$-cut problem is not just {\bf NP}-hard, but also {\bf APX}-hard, i.e., there is a constant $\delta >1$ such that it is {\bf NP}-hard to even approximate the solution to within a ratio of less than $\delta$ to the optimal cost even when restricted to instances with three terminals $(k=3)$ and unit edge costs. The special case of the problem on planar graphs is also {\bf NP}-hard if $k$ is arbitrarily large, but can be solved in polynomial time for every fixed $k$ \cite{marx2012tight}. The complexity of multiway $k$-cut problem arises from the {\em combinatorial} number of ways in which the vertex set $V$ can be partitioned into $k$ feasible sets. 

In many application areas such as social networks, communication networks and epidemic spread networks, the corresponding graph models are time-varying in nature \cite{casteigts2012time, thomas2016review}. For instance, a company wants to maintain its prominence over a certain user base represented by a graph $G(V,E)$. The set of vertices $V$ includes all the existing and potential costumers and the set $E$ of edges denotes the influence of one customer over the other. Such a graph is temporal in nature since the influence of two nodes over one another changes with time depending on several factors such as the number of text messages exchanged in a day. The company identifies $k$ agents $\{s_1,\hdots,s_k\}$ in the set $V$ of vertices and partitions the graph into $k$ subsets such that there is exactly one agent in each subset and the total cut value is minimized. The agents help the company provide services to the customers, as well as, influence them in company's favor by providing customized schemes for each subset. Owing to the dynamic nature of this setting it becomes imperative to find the time-varying minimum multiway $k$-cuts and also, if possible, influence consumer interactions such that cut value remains small (or decreases) with time. In this context the objective is to design control input that externally influence the edge-weight dynamics and determine the minimum multiway $k$-cut at each time instant.

In this paper, we first present an innovative relaxation heuristic for the minimum multiway $k$-cut problem on static digraphs (weighted and directed graphs). The heuristic is based on maximum-entropy-principle (MEP), which in turn, has close analogy to minimum free-energy principle in statistical physics \cite{jaynes1957information}. Using this heuristic we introduce probability distributions on the space of associations between vertices and the terminal points; this constitutes the {\em soft}-partitioning of the graph. We then seek the fairest distribution (that is with maximum entropy) that guarantees the expected cut values to lie below a prescribed upper bound. A sequence of such problems are solved as the values of the upper bound are successively decreased. At the end of these iterations, the resulting distributions are concentrated about the mean values, and these mean values approximate the solutions to the original unrelaxed problem. In the context of data clustering this MEP based heuristic has resulted in the deterministic annealing (DA) algorithm \cite{rose1998deterministic}. DA has been successfully applied to a large class of optimization problems such as, pattern classification \cite{gehler2007deterministic}, image segmentation \cite{mitra2000adaptive}, graph aggregation \cite{xu2014aggregation}, robust speech recognition \cite{rao1997design}, multiple travelling salesman \cite{baranwal2017multiple} and coverage control \cite{xu2014clustering}. However, DA is primarily a resource allocation algorithm and requires the precise knowledge of coordinates of data points that need to be clustered. An advantage of our proposed methodology is that it does not require the knowledge of coordinates of data points provided that a {\em similarity measure} (such as edge-weight matrix in the context of a graph)  between the data points is known. 

In this work, we propose an MEP based algorithm that does away with this limitation of the DA algorithm and is applicable to optimization problems on graphs. In particular, we develop an iterative  algorithm for the minimum multiway $k$-cut problem on {\em weighted} and {\em directed} graphs and present its convergence analysis by exploiting the algebraic structure of the relaxed cost function and non-negativity of the Kullback-Leibler divergence. An important feature of our algorithm is that it is independent of the initialization and is designed to avoid poor local minima. In fact, the algorithm always starts with uniform distribution for the target probability distributions which evolve into appropriate concentrated distributions. For a graph with $N$ vertices, the run-time complexity of the proposed algorithm  is  $\sim O(kIN^3)$, where $N$ and $I$ denote the number of vertices and iterations respectively. Empirical evaluations on planar and non-planar graphs suggest that our algorithm outperforms the approximation algorithms described in \cite{dahlhaus1994complexity, cualinescu1998improved, naor20012}. In fact, these approximation algorithms are known to result in highly suboptimal solutions on specific graph instances. We show that our approach returns optimal solutions on such otherwise challenging instances. We further demonstrate that our algorithm is also capable of handling very large graph instances (comprising of image pixels as nodes) for interactive foreground-background segmentation.

In the context of determining the minimum multiway $k$-cut on a dynamic graph there could be several simple approaches; one of the straightforward method is the {\em frame-by-frame} approach where we solve for the minimum cut at each time instant. However such a methodology is computationally expensive and non-viable for very large graphs. In our proposed methodology we first choose an energy function that captures the multiway $k$-cut value as well the external control effort. We then treat this function as a control Lyapunov function to design the control input, which determines the manipulable edge-weight dynamics and the time-varying minimum multiway $k$-cuts. We show that this approach is not conservative, that is, under the assumption of feasibility, it always results in a control law that ensures the non-positiveness of the time derivative of the energy function. Our simulations demonstrate reduction in computational times by over $25$ times in comparison to the frame-by-frame approach.

\section{Problem formulation in Static Graphs}\label{sec:Problem_Def}
For a given weighted directed graph $G=(V,E,W)$ with vertex set $V, |V| = N\in\mathbb{N}$, set $E\subseteq V\times V$ of edges, edge weight matrix $W=[w_{lm}]\in \mathbb{R}_{\geq 0}^{N\times N}$ with $w_{lm}$ as the weight of the edge from vertex $l$ to $m$ and $w_{lm} = 0$ $\forall$ $(l,m)\notin E$, and a set of terminals $S = \{s_1,s_2,\dots,s_k\}\subseteq V$, the minimum multiway $k$-cut problem is defined as
\begin{align}\label{eq:problem}
	&\min_{\substack{\{A_1,A_2,\hdots,A_k\},\\ A_j\subset V}}\frac{1}{2}\sum_{j_1=1}^{k}\sum_{\substack{j_2=1\\j_2\neq j_1}}^{k}\sum_{\substack{l\in{A_{j_1}}\\m\in{A_{j_2}}}}w_{lm} \nonumber \\
	&\text{s.t.}\quad s_1\in A_1,\hdots,s_k\in A_k, \quad  \cup_{i=1}^kA_k=V,\nonumber \\
	&\text{and}\qquad A_j\cap A_l=\phi, \ \text{for all}\  j\neq l.
\end{align}

Note that, a partition $\{A_1,\dots,A_k\}$ of $V$ results in $k$ disconnected components (subgraphs) of graph $G$ (see Fig. \ref{fig:multiwayCut}a). A component $A_j$ contains all the vertices $i$ in $V$ that are in the same subgraph as the terminal $s_j$. Here a {\em multiway $k$ cut} represents the set of edges whose vertices belong to distinct components (i.e. distinct $A_j$'s). For instance the set $\{e_1,e_2,e_3\}$ is a multiway $k$ cut for the graph in Fig. \ref{fig:multiwayCut}(a) where set of terminals is $\{s_1,s_2,s_3\}$. The above optimization problem seeks a partition that minimizes the total {\em cut weight}, i.e. the cumulative weight of all the edges in the multiway $k$ cut.

We reformulate the optimization problem in (\ref{eq:problem}) by introducing $N\times k$ {\em soft} decision variables $\{p(j|i)\}$. These variables describe soft partitions of the graph. Here $p(j|i)\in[0,1]$ denotes {\em a probabilistic} association of the vertex $i$ with terminal $s_j$. We also require that for each vertex $i$, $\sum_j p(j|i)=1$, i.e., $\{p(j|i)\}$ ascribes a probability distribution over all feasible associations over vertex set $V$. In the case of hard partitions in (\ref{eq:problem}) all the decision variables $\{p(j|i)\}$ will be taking values either $0$ or $1$. In this case a component $A_j$ corresponds to a set $\{i\in V: p(j|i)=1\}$. A relaxation of the optimization problem (\ref{eq:problem}) is given by 
\begin{align}\label{eq:D_multiway_cut}
	\min_{\substack{\{p(j|i)\}\\1\leq i\leq N\\1\leq j\leq k}} &D \triangleq \frac{1}{2}\sum_{\substack{j_1,j_2=1\\j_2\neq j_1}}^{k}\sum_{l,m=1}^N p(j_1|l)p(j_2|m)w_{lm} \nonumber\\
	\text{s.t.}\quad &p(j_1|s_1)=\dots=p(j_k|s_k)=1,
\end{align}
where we set  $W \leftarrow W+\lambda\mathbb{I}_N$. The inclusion of a constant parameter $\lambda$ or equivalently the {\em regularizer} term in $D$ is explained as follows.   Note that in formulation (\ref{eq:problem}), self-loop edges with node  weights  $w(l,l)$ can not be in a multiway $k$-cut, and therefore the solution to this problem is {\em independent} of these weights, that is, independent of the diagonal entries of the matrix $W$. However in the proposed  relaxation, since each vertex has partial membership with respect to different terminals $s_j$, the self-loop edges can become a part of the multiway $k$-cut; so we replace the edge-weight matrix $W$ by $W+\lambda I_N$, i.e.
\begin{align}\label{eq:LambdaW}
    W\longleftarrow W+\lambda I_{N},
\end{align}
where $I_N$ is a $N\times N$ identity matrix and $\lambda>0$ is large enough to make sure that self-loop edges are not included in the multiway $k$-cut; the choice of $\lambda$ is discussed in Section \ref{sec:proof_k_2}.

\begin{remark}
Note that without loss of generality we can assume the edge-weight matrix $W$ to be symmetric. In fact, solution to the optimization problem in (\ref{eq:D_multiway_cut}) for a given edge-weight matrix $W$ is identical to solution of a similar optimization problem with edge-weight matrix $\tilde{W}:=0.5(W+W^T)$ which can be verified by substituting $\sum_{j_2\neq j_1}p(j_2|m)=1-p(j_1|m)$ in expression of $D$ in (\ref{eq:D_multiway_cut}).
\end{remark}
\section{Solution to Multiway $k$-cut in Static Graphs}\label{sec:Problem_Soln}
In the proposed approach instead of directly solving (\ref{eq:D_multiway_cut}) we use the Maximum Entropy Principle (MEP) \cite{jaynes1957information} to determine the distributions $\{p(j|i)\}$ that ensure the relaxed cost function $D$ in (\ref{eq:D_multiway_cut}) is less than or equal to a constant $d_0>0$. More specifically, this principle states that of all the probability distributions that satisfy a given set of constraints on expected values of functions of a random variable, choose the one that maximizes the Shannon entropy $H\left(\{p(j|i)\}\right)$. Accordingly in our case, the MEP would solve $\max ~H\left(\{p(j|i)\}\right)$ under the constraint that $D\leq d_0$, where $D$ is given in (\ref{eq:D_multiway_cut}) and the Shannon entropy term is given by
\begin{align}\label{eq:H}
\begin{small}
H \triangleq -\sum\limits_{i=1}^N\sum\limits_{j=1}^k p(j|i)\log{p(j|i)}.
\end{small}
\end{align}
The equivalent Lagrangian is thus defined as
\begin{align}\label{eq:L}
\begin{small}
L\triangleq D-d_0 - \frac{1}{\beta}H,
\end{small}
\end{align}
where Lagrange multiplier $\beta$ controls the trade-off between minimizing cost function $D$ and maximizing entropy $H$.

In the expression of the Lagrangian (\ref{eq:L}) we refer to the Lagrange multiplier $1/\beta$ as temperature and L as {\em free energy} because of their close analogies to statistical physics (where free energy is enthalpy ($D$) minus the temperature times entropy ($TH$)). Note that the free energy $L$ can also be viewed as a relaxation of the cost function (\ref{eq:D_multiway_cut}). In fact, as $\beta~\to~\infty$, we note that $L\rightarrow D$. Using the fact that $\sum_{j_2\neq j_1}^{k}p(j_2|i)=1-p(j_1|i)$ and $d_0$ is a constant, the effective Lagrangian in (\ref{eq:L}) is given by
\begin{align}\label{eq:F}
L &= \frac{1}{2}\sum_{l,m=1}^Nw_{lm}\Big(1-\sum_{j=1}^k p(j|l)p(j|m)\Big)\nonumber\\
&\quad + \frac{1}{\beta} \sum\limits_{i=1}^N \sum\limits_{j=1}^k  p(j|i)) \log{p(j|i)}.
\end{align}

The above reformulation is critical to the convergence analysis described in Section \ref{sec:proof_k_2}.  By setting $\frac{\partial L}{\partial p(j|i)}=0$ toward minimizing (local) $L$ with respect to $p(j|i)$ yields
\begin{align}\label{eq:P}
\begin{small}
p(j|i) = \frac{\exp\big\{ \beta \big(\sum\limits_{m=1}^{N}p(j|m) w_{im}\big)\big\}}{Z_i},
\end{small}
\end{align}
where the normalizing constant $Z_i$ is given as
\begin{align}\label{eq:Normalize}
\begin{small}
Z_i = \sum_j\exp\big\{ \beta \big(\sum\limits_{m=1}^{N}p(j|m) w_{im}\big)\big\}.
\end{small}
\end{align}
We then substitute the expression (\ref{eq:P}) of $p(j|i)$ in the free-energy $L$ (\ref{eq:F}) to obtain

\begin{align}\label{eq:F_final}
L & \triangleq  \frac{1}{2}\sum_{l,m=1}^N w_{lm}\Big(1-\sum_{j=1}^k p(j|l)p(j|m)\Big)\nonumber\\
&\quad-\frac{1}{\beta}\sum_{i=1}^N \log\Big(\sum_{j=1}^k \exp{\big(\sum_{m=1}^N p(j|m)w_{im}\big)}\Big).
\end{align}

The essence of the MEP-based approach lies in successive evaluations of Gibbs distribution in (\ref{eq:P}). Note that from (\ref{eq:F_final}) minimizing $L$ at small values of $\beta$ is equivalent to maximizing entropy $H$, which in turn corresponds to uniform distribution. As $\beta$ is gradually increased, minimization of Lagrangian in (\ref{eq:L}) puts more weight on minimization of the cost function (\ref{eq:D_multiway_cut}) and as evident from (\ref{eq:P}) results in {\em hard} (0-1) associations $\{p(j|i)\}$ as $\beta\rightarrow\infty$. This process of gradual {\em cooling} is referred to as {\em annealing} in the statistical physics literature. Also observe that increasing the Lagrange parameter $\beta$ is equivalent to solving the same MEP problem with a decreased value of $d_0$ \cite{jaynes1957information}. Thus as $\beta\to\infty$, the algorithm seeks the minimum value  for the cost $D$ in (\ref{eq:D_multiway_cut}). In the proposed algorithm (see Algorithm \ref{alg:alg1}) we minimize the free-energy $L$ through fixed-point iterations in (\ref{eq:P}) at successively increasing values of the annealing parameter $\beta$. Note that at each value of the annealing parameter $\beta$, the algorithm executes the following two steps to solve the equation (\ref{eq:P})
\begin{align}\label{eq:steps}
\text{{\bf Step 1}:} \qquad &\sigma_{lj}^+\leftarrow p(j|l), \nonumber \\
\text{{\bf Step 2}:} \qquad &p^+(j|l) \leftarrow \frac{\exp\left\{ \beta \left(\sum\limits_{m=1}^{N}\sigma^+_{mj} w_{lm}\right)\right\}}{\sum_j\exp\left\{ \beta \left(\sum\limits_{m=1}^{N}\sigma^+_{mj} w_{lm}\right)\right\}}, \nonumber \\
&\forall l\in\{1,\dots,N\}\setminus S, j\in\{1,\dots,k\}
\end{align}
\begin{algorithm}
	\caption{Algorithm for minimum MultiwayCut}
	\begin{algorithmic}
		\renewcommand{\algorithmicrequire}{\textbf{Input:}}
		\renewcommand{\algorithmicensure}{\textbf{Output:}}
		\REQUIRE $G = (V,E)$, $w:E\to\mathbb{R}^+$, $S=\{s_1,\dots,s_k\}$
		\ENSURE $\{p(j|i)\}$
		\\ \textit{Initialization}:
		\\$p(j|j)\gets 1\quad\forall s_j\in S$,
		\quad$p(j|i)\gets \frac{1}{k} \quad \forall i\not\in S, \forall j$
		\\$\beta \gets \beta_{\min}$
		\\ \textit{Annealing Process}
		\WHILE {$\beta < \beta_{\max}$}
		\STATE \textit{Fixed-Point Iterations}
		\WHILE {until convergence}
		\STATE Update $\{p(j|i)\}$ as in (\ref{eq:P}) $\forall i\not\in S$
		\ENDWHILE
		\STATE Increment $\beta$
		\ENDWHILE
		\RETURN $\{p(j|i)\}$ 
	\end{algorithmic}
	\label{alg:alg1}
\end{algorithm}

{\em Time-complexity of the proposed algorithm}: The main complexity of this algorithm stems from the matrix multiplication in the fixed point iteration scheme. For a graph $G = (V,E)$ with $|V|=N$, there is a total $Nk$ association probability parameters $\{p(j|i)\}$, that need to be estimated at each $\beta$ iteration. Note that the {\em batch} update equation in (\ref{eq:P}) requires multiplying the two matrices $W$ and $[p(j|i)]_{N\times k}$. This multiplication operation runs in $O(N|E|k)$ time (total of $N|E|$ operations for each partition associated with multiplying non-zero elements of edge-weight matrix). Thus the run-time complexity of the proposed algorithm is $\approx O(kIN^3)$ where $I$ accounts for the number of $\beta$ iterations and the fixed-point iterations. 

\begin{figure}
    \centering
    \includegraphics[width=0.6\columnwidth,height=0.5\columnwidth]{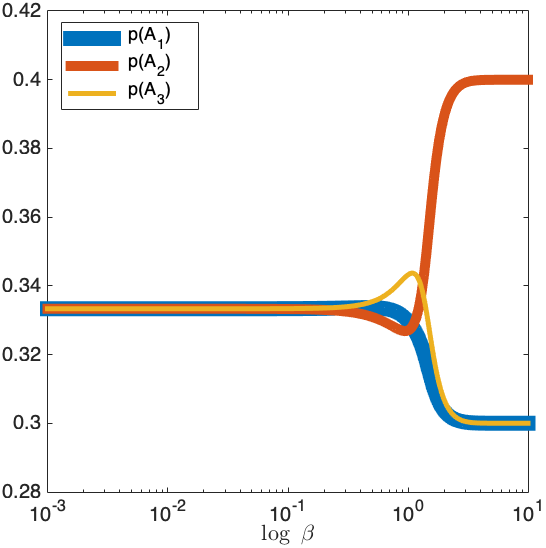}
    \caption{{\small Phase Transition in Algorithm 1 applied to problem stated in Fig. \ref{fig:multiwayCut}(b) with unit edge-weights.}}
    \label{fig:PT}
\end{figure}

{\em Remark : }Usually the MEP-based heuristics developed for solving the combinatorial optimization problems undergo several {\em phase transitions} \cite{rose1998deterministic} as the parameter $\beta$ increases from zero to a large number. In the case of the Algorithm \ref{alg:alg1}, these phase transitions that occur at certain critical $\beta_{cr}$'s, correspond to an abrupt change in the weight of the partitions ${A_j}$ defined as $p(A_j) := \sum_i p(j|i)/N$. Empirical evaluations suggest that ${p(A_j)}$ does not change much between two consecutive phase transitions. Hence in our simulations we anneal the parameter $\beta$ geometrically, i.e. $\beta_t = \gamma^{t}\beta_0$ where $\gamma>1$ and the effective number of $\beta$ iterations is $\log\frac{\beta_{\max}}{\beta_{\min}}$ which is small. In fact in all our simulations $\log\frac{\beta_{\max}}{\beta_{\min}}\lessapprox 10$. Fig. (\ref{fig:PT}) demonstrates the phase transition phenomenon observed in our simulation when Algorithm \ref{alg:alg1} is applied to the problem stated in Fig. \ref{fig:multiwayCut}(b) with unit edge-weights. In our ongoing work we are solving for analytical conditions that determine the critical $\beta_{cr}$ at which the phase transition phenomenon occurs.

\section{Convergence analysis}\label{sec:proof_k_2}
In this section we provide a proof that the Lagrangian $L$ in (\ref{eq:F}) converges to a local minimum under the two-step iterations specified by (\ref{eq:steps}). More specifically we show that for a fixed value of the Lagrange parameter $\beta$, every successive iteration of (\ref{eq:steps}) decreases the effective Lagrangian $L$ and since $L$ is lower bounded it converges to a local minimum. We use $P\in\mathbb{R}^{N\times k}$ to denote the matrix of associations $[p(j|i)]$.

\textit{Claim : The Lagrangian $L(P)$ in (\ref{eq:F}) converges to a local minimum under the fixed point iterations (\ref{eq:P}) (equivalently the two step iterations in \ref{eq:steps}) in Algorithm \ref{alg:alg1}.}

\textit{Proof}. We first show that $L(P)$ is decreasing under the two-step iterations (\ref{eq:steps}); that is we show that $\Delta(P,P^+):=L(P)-L(P^+)>0$ where $P^+$ is obtained after successive executions of steps 1 and 2 in (\ref{eq:steps}). Towards this end, we first construct a function $\Gamma(\zeta,\eta):\mathbb{R}^{N\times k}\times \mathbb{R}^{N\times k}\rightarrow \mathbb{R}$ as

{
\small
\begin{align}\label{eq:E}
\Gamma(\zeta,\eta) \triangleq \frac{1}{2}\sum\limits_{l,m,j=1}^{N,N,k}\Big(\eta_{lj}-2\zeta_{lj}\Big)\eta_{mj} w_{lm} + \frac{1}{\beta} \sum\limits_{i,j=1}^{N,k} \zeta_{ij}\log{\zeta_{ij}};
\end{align}}

here observe that $\Gamma(P,P)=L(P)$. Therefore 
\begin{align}
\begin{small}
    \Delta(P,P^+) =\Gamma(P,P)-\Gamma(P^+,P^+),
\end{small}
\end{align}
using the fact that $\Gamma(P,P)=L(P)$. To show that $\Delta(P,P^+)\geq 0$ we re-write it as
\begin{align}
\begin{small}
	\Delta(P,P^+) =  \underbrace{\Gamma(P^+,P)-\Gamma(P^+,P^+)}_{\Delta_1} + \underbrace{\Gamma(P,P)-\Gamma(P^+,P)}_{\Delta_2},
\end{small}
\end{align}
and then show that both $\Delta_1$ and $\Delta_2$ are non-negative, i.e. $\Delta_1\geq 0\text{ and }\Delta_2\geq 0$. 
We show that $\Delta_1\geq 0$ by showing that $\Gamma(P^+,\sigma)$ (a quadratic function of $\sigma$ with fixed $P^+$) achieves its minimum when $\sigma=P^+$.  The stationary point is obtained by   
setting $\frac{\partial \Gamma}{\partial \sigma_{ij}}\big|_{(P^+,\sigma)}=0$ in (\ref{eq:E}), which yields
\begin{align}\label{eq:E_der}
\sigma_{ij} &= p^+(j|i). 
\end{align}
This solution $\sigma=P^+$ is a minimum if the Hessian $\frac{\partial^2\Gamma(P^+,\sigma)}{\partial\sigma^2}\big|_{\sigma=P^+}=W$ is positive-definite. Choosing $\lambda$ in (\ref{eq:LambdaW}) such that $\lambda =\sum_i \sum_{j\neq i} w_{ij}$, ensures positive-definiteness of the $W$ matrix by the Gershgorin circle theorem \cite{gershgorin1931uber}. Since $P^+$ is a minimizer, setting $\sigma=P\neq P^+$ results in $\Delta_1\geq 0$.

To show that $\Delta_2\geq 0$, note that from (\ref{eq:E}) we have

{\small
\begin{align}\label{eq:v1}
&\Delta_2=-\frac{1}{2}\sum \limits_{lmj}\left(p(j|l)-p^+(j|l)\right)p(j|m) w_{lm}\nonumber\\
&+\frac{1}{\beta}\sum_{lj}\big(p(j|l)\log p(j|l)- p^+(j|l)\log p^+(j|l)\big).
\end{align}}

Also from  taking logarithm on both sides of (\ref{eq:P}), we have
\begin{align}\label{eq:v2}
\begin{small}
\sum \limits_m p(j|m) w_{lm}=\frac{1}{\beta}\log Z_l+\frac{1}{\beta}\log p(j|l),
\end{small}
\end{align}
\noindent where $Z_l$ is given by (\ref{eq:Normalize}). Using (\ref{eq:v2}) and the fact that $\sum_j p(j|l)=1$ $\forall$ $l\in V$ we obtain
\begin{align}\label{eq:v3}
&\sum \limits_{lj}\left(p(j|l)-p^+(j|l)\right)\sum\limits_m p(j|m) w_{lm}\nonumber\\
&=\frac{1}{\beta}\sum \limits_{lj}\left(p(j|l)-p^+(j|l)\right)\log p^+(j|l).
\end{align}
On substituting (\ref{eq:v3}) in (\ref{eq:v1}), we obtain
\begin{small}
\begin{align}
\Delta_2 = \frac{1}{\beta}\sum\limits_{l,j=1}^{N,k}p(j|l)\log{\frac{p(j|l)}{p^+ (j|l)}}
=\frac{1}{\beta}\sum\limits_{l=1}^{N}D_{KL}(P_l||P^+_l)\geq 0,
\end{align}
\end{small}
\noindent where $P_l$ and $P^+_l$ are the $l$th rows of the matrices $P$ and $P^+$ respectively, and $D_{KL}(P_l||P^+_l)$ represents the Kullback-Leibler measure. Consequently, $\Delta(P,P^+)\geq0$. Thus we have shown that $L(P)$ decreases as a result of the two-step iteration in (\ref{eq:steps}), and since $L(P)$ is bounded from below, the fixed point iterations in (\ref{eq:P}) converge to a local minimum of $L(P)$.

\section{Extension to Dynamic Graphs}\label{sec:DynamicCase}
Here we consider a dynamic digraph $G(t)=(V,E,W_0,\Sigma)$ where $V$, $|V|=N\in\mathbb{N}$ denotes the set of vertices, $E\subseteq V\times V$ denotes the set of edges, $W_0=[w_{0_{lm}}]\in \mathbb{R}^{N\times N}$ denotes the edge-weight matrix with $w_{0_{lm}}=0$ $\forall$ $(l,m)\notin E$, and $\Sigma$ denotes the dynamics of the time-varying edge-weight matrix $W(t)=[w_{lm}(t)]$ given by $\dot{W} = \mathcal{F}(W,U)$ where $W(t_0)=W_0$ and $\dot{w}_{lm}(t)=0$ $\forall$ $(l,m)\notin E$. We assume that the matrix function $\mathcal{F}:=[f_{ij}]$ is known a priori and  defined as $\small{f_{ij}:\mathbb{R}^{N\times N}\times \mathbb{R}^{N\times N}\to \mathbb{R}}$ where $f_{ij}$ belong to the class of continuously differentiable $\small{C^1(\mathbb{R}^{N\times N})}$ functions. Here $U(t)=[u_{ij}(t)]\in\mathbb{R}^{N\times N}$ denotes the external control input  and is defined as $u_{ij}:\mathbb{R}\rightarrow \mathbb{R}$. The function $u_{ij}$ belongs to the class of continuously differentiable $C^1(\mathbb{R})$ functions.
\begin{remark}\label{rem:DynStr}
The matrix function $\mathcal{F}(W,U)$ encodes the information whether a particular edge $(l,m)\in E$ in the graph is manipulable or not. 
\end{remark}

The objective here is two-fold (a) design a control $U(t)$ such that a modified multiway $k$-cut value in (\ref{eq:DynamicP}) is minimized at every time instant $t$, and (b) determine the time-varying multiway $k$-cut for the dynamically evolving graph $G(t)$ given the set of terminal $S=\{s_1,\hdots,s_k\}\subseteq V$. The corresponding dynamic optimization problem is given by
\begin{align}\label{eq:DynamicP}
    \min_{\substack{\{A_j\}, A_j\subset V\\\{u_{ij}(t)\in C^1\}}}\quad &\frac{1}{2}\sum_{j_1=1}^k\sum_{\substack{j_2=1\\j_2\neq j_1}}^k\sum_{\substack{l\in A_{j_1}\\m\in A_{j_2}}}\big(w_{lm}(t)\big)^2 + \mu \|U(t)\|^2_{F}\nonumber\\
    \text{s.t. }\quad & s_j \in A_j \forall~j\in\{1,\hdots,k\}, \cup_{j=1}^k A_j=V,\nonumber\\
    & A_j\cap A_l = \phi \forall ~j\neq l,\nonumber\\
    & \dot{W}(t) = \mathcal{F}(W,U), \quad W(t_0) = W_0,
\end{align}
where $\mu\in \mathbb{R}$ is a user-defined parameter, and $\|\cdot\|_{F}$ is the Frobenius norm. The first part of this objective function corresponds to the multiway $k$-cut cost, where each edge-weight is replaced by its square, and the second part represents a penalty on the control effort. The user-defined parameter $\mu$ regulates the relative weight given to the control effort $\mu\|U(t)\|_\mathrm{F}^2$. Observe that the edge-weight dynamics described by $\dot{W}(t)=\mathcal{F}(W,U)$ also allows for edge-weights to possibly become negative. This issue is addressed by modifying the cost function to have squared weight terms $\{w_{lm}^2\}$ instead of $\{w_{lm}\}$. In this way, the cost function penalizes only the magnitude (square) of the edge weights regardless of their signs
and the resulting cuts include only those edges whose
weights have small magnitude; a feature required in most applications. We discuss the case with the original (not modified) cost function in remark \ref{rem:PutWBack}. The dynamic optimization problem in (\ref{eq:DynamicP}) inherits the computational complexity of the static-problem; which is further worsened by the dynamical aspect of the problem. Solving for $U(t)$ at each time-instant while ensuring that the function $U(t)$ is smooth is one of the main contributors to the additional complexity. 

One straightforward method is the {\em frame-by-frame} approach which disregards this constraint. Here we set $U(t)=0$ for all $t$ and solve the minimum multiway $k$-cut problem using the Algorithm \ref{alg:alg1} (for the cost function in (\ref{eq:DynamicP})) at every time instant $t$. A disadvantage of this approach is that, if the time interval $\Delta t_p$ between two successive runs of the Algorithm \ref{alg:alg1} is short, then the overall approach is computationally expensive. Also the frame-by-frame analysis does not exploit the information available from the previous time instances to determine the minimum multiway $k$-cut at the current time instant. On the other hand, if the time interval $\Delta t_p$ between two successive instances is large, then the algorithm cannot account for dynamics in this time interval and the resulting cut may be correspondingly large. 

We propose an alternative method, where instead of solving directly the dynamic optimization problem in (\ref{eq:DynamicP}), we address its objectives of minimizing the multiway $k$-cut and the control effort. Here we consider an energy-like function $F(W(t),U(t))$ and design a dynamic control law $\dot{U}(t)=\mathcal{V}(W,U)$ such that $\dot{F}(W,U)\leq 0$. The function is given by 
\begin{align}\label{eq:Lyapunov}
{\small
F(W(t),U(t)):= L_1(W(t)) + \mu \|U(t)\|_\mathrm{F}^2,}
\end{align}
where $L_1(W(t))=$
{\small
\begin{align}\label{eq:L1}
&\frac{1}{2}\sum_{l,m=1}^N w_{lm}^2\big(1-\sum_{j=1}^k p(j|l)p(j|m)\big) + \frac{1}{\beta}\sum_{i,j=1}^{N,k}p(j|i)\log p(j|i)
\end{align}}
\noindent is the effective Lagrangian (\ref{eq:F}) modified with the squared edge-weights. Since this Lagrangian is a close approximation of the multiway $k$-cut value described in (\ref{eq:DynamicP}) (especially for high values of $\beta$ as seen in Section \ref{sec:Problem_Soln}); its inclusion in $F(W,U)$ addresses the objective of having {\em small} multiway $k$-cut values. 

The time-derivative $\dot{F}(W,U)$ is given by 
\begin{align}\label{eq:Vdot}
{\small
\dot{F} = e_N^T\Big(\Phi\circ W(t)\circ\mathcal{F}(W,U)(t)+2\mu U(t)\circ\mathcal{V}(t)\Big)e_N,}
\end{align}
where $\dot{U}(t) = \mathcal{V}(W,U)$ represents a dynamic control law, $e_N\in \mathbb{R}^{N}$ is a column of all $1$'s, $\Phi=(\phi_{lm})\in\mathbb{R}^{N\times N}$ such that $\phi_{lm}=1-\sum_{j=1}^k p(j|l)p(j|m)$ and $\circ$ denotes the Hadamard product. We exploit the affine dependence of $\dot{F}$ on $\mathcal{V}$ in (\ref{eq:Vdot}) to make $\dot{F}$ non-positive analogous to control based on control Lyapunov functions \cite{sepulchre2012constructive, sontag1983lyapunov, sontag1989universal}. Specifically we choose 
\begin{align}\label{eq:cont_dyn}
{\small
    \mathcal{V} = -\Bigg[C_0 + \frac{\alpha+\sqrt{\alpha^2+(4\mu \|U(t)\|_\mathrm{F}^2)^2}}{4\mu \|U(t)\|_\mathrm{F}^2}\Bigg]U(t),}
\end{align}
where $\alpha = 2e_N^T\Phi\circ W\circ \mathcal{F}(W,U)e_N$ and $C_0>0$. The following theorem summarizes the consequences of this design:
\begin{theorem} 
\[\]\vspace{-0.40in}
\begin{enumerate}
\item $F(W,U)$ is lower bounded; more specifically $F(W,U)+\dfrac{1}{\beta}N\log k\geq 0$. 
\item If there exists a dynamic controller $\dot U=\bar{\mathcal{V}}(W,U)$ such that 
\begin{enumerate}
\item $\mathcal{V}(W,U)$ is locally Lipschitz,
\item $\dot{F} = \frac{1}{2}\alpha + 2\mu U(t)\circ\bar{\mathcal{V}}\leq 0$, where $\alpha = 2e_N^T\Phi\circ W\circ\mathcal{F}(W,U)e_{N}$,
\end{enumerate}
then the control design $\dot U=\mathcal{V}(W,U)$ in (\ref{eq:cont_dyn}) is such that 
\begin{enumerate}
\item $\bar{\mathcal{V}}(W,U)$ is bounded and locally Lipschitz
\item $\dot{F} = \frac{1}{2}\alpha + 2\mu U(t)\circ {\mathcal{V}}\leq 0$, where $\alpha = 2e_N^T\Phi\circ W\circ\mathcal{F}(W,U)e_{N}$
\item $\dot F\rightarrow 0$,   $ U \rightarrow 0$, and  $\dot L_1 \rightarrow 0$ as $t\rightarrow \infty$.
\end{enumerate}
\end{enumerate}
\end{theorem}

Please refer to Appendix in Section \ref{sec:appendix} for the proof of the above theorem.

\begin{remark}
From the above theorem we show that our control design $\dot U=\mathcal{V}$ achieves the objectives of (\ref{eq:DynamicP}), that is $\dot F\leq 0$, only when {\em there exists atleast one} Lipshitz control design $\bar{\mathcal{V}}$ such that $\dot{F}\leq 0$. The conditions under which such a $\bar{\mathcal{V}}$ exists is a difficult problem in itself and is part of our ongoing work. Essentially this theorem shows that our approach is not conservative; it guarantees $\dot F\leq 0$, whenever it is possible (by any control design) under the dynamics $\dot W = \mathcal{F}(W,U)$.  
\end{remark}

\begin{remark}\label{rem:PutWBack}
The proposed control design approach when applied to the case with the original  (not modified) cost function, that is is with $L$ in (\ref{eq:F}); yields similar results. In fact the resulting dynamic control law is Lipshitz and is such  that $\dot F\leq 0$ under the assumption that there exists a $\bar{\mathcal{V}}$ that satisfies the conditions in the above theorem.
\end{remark}
\section{Illustrative Examples}\label{sec:sim_results} 
\begin{figure}
	\centering
	\includegraphics[width=\columnwidth]{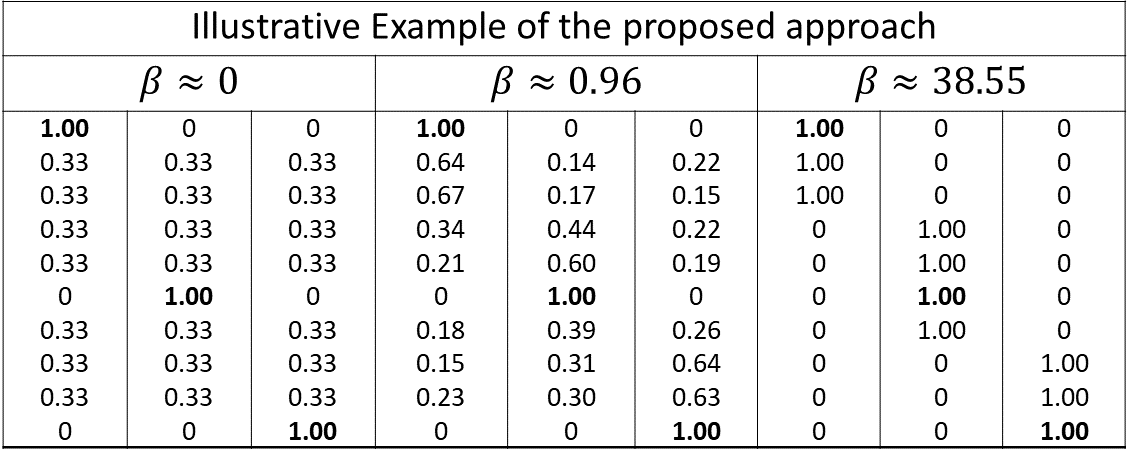}
	\caption{{\small Illustrative example of minimum $3$-cut problem. The columns indicate the association matrices $p(j|i)$ at different $\beta$ values. The `bold' numbers indicate the associations of terminals, i.e., $p(j=1|i=1)=1, p(j=2|i=6)=1$ and $p(j=3|i=10)=1$. As the algorithm progresses, the associations of remaining nodes harden, thereby resulting in optimal cut.}}
	\label{table:tab1}
\end{figure}

We demonstrate our algorithm for the example shown in Fig. \ref{fig:multiwayCut}b. The Algorithm \ref{alg:alg1} determines the partitioning $\{A_j\}_{j=1}^3$ as $A_1=\{1,2,3\}$, $A_2=\{8,9,10\}$, $A_3=\{4,5,6,7\}$ and the resulting set of cut-edges as $\{(2,9),(3,4),(7,8)\}$. It is easy to verify that the above solution is optimal. Fig. \ref{table:tab1} shows a typical run of our algorithm at increasing $\beta$ values.
\begin{figure*}
    \centering
    \includegraphics[width=1.9\columnwidth]{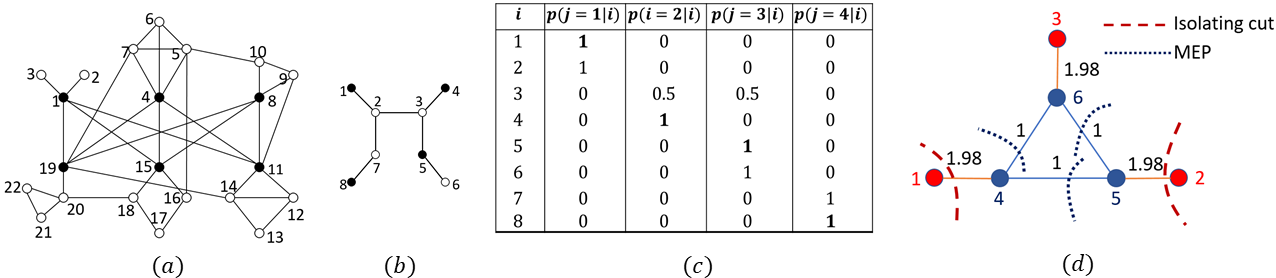}
    \caption{{\small (a) Example of a $22$ nodes non-planar graph with $6$ terminals, $S = \{1,4,8,11,15,19\}$. The proposed algorithm results in an optimal cut. (b) An $8$ node graph with multiple permissible optimal $4$-cuts. (c) Final association matrix for the $4$-cut problem shown in (b) part. (d) Toy-example with $2k=6$ vertices where $k$ of the vertices form a cycle where each edge weight is equal to 1, and each other vertex is connected to exactly one vertex on the cycle with an edge weight of 1.98. The isolating cut heuristic \cite{dahlhaus1994complexity} results in a cut value of 3.96, where as our MEP-based algorithm results in optimal cut.}}\label{fig:multiwayCutEx}
\end{figure*}
\begin{figure}
    \centering
    \includegraphics[width=0.6\columnwidth]{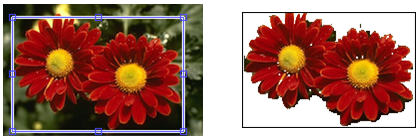}
    \caption{{\small Illustration of the proposed multiway $k$-cut approach to background-foreground segmentation. Here each pixel in the bounding box represents a node. The terminal nodes (corresponding to foreground and background) are interactively provided by the user once, and the algorithm successively evaluates minimum $st$-cut to refine foreground-background segmentation. Edge-weights are obtained as functions of pixels' RGB values.}}
    \label{fig:grabcut}
\end{figure}

As seen in Fig. \ref{table:tab1}, at very low values of $\beta$, the association probability $p(j|i)\equiv 1/3$ for every $i\notin S$, which in turn corresponds to maximizing randomness (Shannon entropy) of the solution. However as the randomness is gradually decreased by increasing $\beta$, the probabilities start becoming {\em non-uniform} and exhibit preferential association to a specific terminal $s_j$. In the limiting case, i.e., at large values of $\beta$, the algorithm results in hardened probabilities ($0-1$ associations). Thus an optimal cut is obtained. Note that in this example, $\beta$ is increased {\em geometrically} from $0.01$ to $40$, i.e., the algorithm provides for very fast $\beta$ scheduling. 

{\em Minimum multiway $k$ cut for a non-planar graph}: We now consider an instance of a $22$-node {\em non planar} graph, shown in Fig. \ref{fig:multiwayCutEx}a with unit edge-weights. In this example we consider a $6$-cut problem, whose set of terminals is specified as $S=\{1,4,8,11,15,19\}$. Our MEP algorithm results in a partition of the underlying graph, given by $\{1,2,3\}$, $\{4,5,6,7\}$, $\{8,9,10\}$, $\{11,12,13,14\}$, $\{15,16,17,18\}$ and $\{19,20,21,22\}$ with a cut value of 15; which is indeed optimal and can be easily verified. On the other hand, the isolating cut heuristic \cite{dahlhaus1994complexity} for this randomly generated instance results in a cut solution with a value of 16. A similar observation is made on other randomly generated instances, where our algorithm results in optimal cut values (whenever verifiable). Moreover, the total run-time for the example in Fig. \ref{fig:multiwayCutEx}a for a naive implementation of the proposed approach in MATLAB is $<1$s on an Intel i7-4790 CPU @ 3.60 GHz.

{\em Non-unique optimal cuts}: In both the examples described above, the resulting optimal cuts are indeed {\em unique}. We now therefore consider a scenario with more than one permissible optimal cuts. Our algorithm identifies the multiplicity of optimal cuts, which is reflected in the final association matrix $\{p(j|i)\}$. Fig. \ref{fig:multiwayCutEx}b shows an example of a $8$-node graph with unit edge-weights and multiple permissible optimal $4$-cuts . Executing the proposed MEP-based algorithm results in the following association matrix. As shown in Fig. \ref{fig:multiwayCutEx}(c), node $3$ can be included in either $A_2$ or $A_3$ without affecting the value of the optimal cut, and therefore the association probabilities $p(j=2|i=3) = p(j=3|i=3) = 0.5$. Assigning node $3$ to any of the two partitions still results in a feasible, yet optimal cut.

{\em Comparison on challenging examples}: Fig. \ref{fig:multiwayCutEx}d shows the performance of our MEP based approach on a toy-example (which can be generalized to any number of nodes. Here heuristic such as isolating cut \cite{dahlhaus1994complexity}, fails to identify an optimal cut for the simplest such scenario. In fact, it has been analytically shown that the isolating cut heuristic results can not give a solution within $2(1-1/k)$ of the optimal solution \cite{dahlhaus1994complexity}. However, our algorithm finds the optimal multiway $k$-cut; thus achieving a cut which is impossible to obtain using the isolating cut method.

{\em On Large Graphs: }We also test our approach on large graphs with the number
of nodes as large as $\sim 25000$ (corresponding to the size
of the bounding box 150$\times$160 pixels). Figure 4 shows the
results of our implementation of the interactive foreground-background segmentation (GrabCut \cite{rother2004grabcut}) using the proposed MEP approach. In GrabCut, an image is represented as a graph of pixels (nodes), where edge-weights capture differences in intensities between neighboring pixels. Users are required to demarcate the approximate foreground region using a bounding box as shown in Fig. \ref{fig:grabcut}. The proposed approach is then employed to segment and refine the foreground through successive minimization of {\em s-t} cuts of the resulting graphs. The implementation results in effective segmentation of foreground and background

\begin{figure}
    \centering
    \includegraphics[width=\columnwidth,height=0.30\columnwidth]{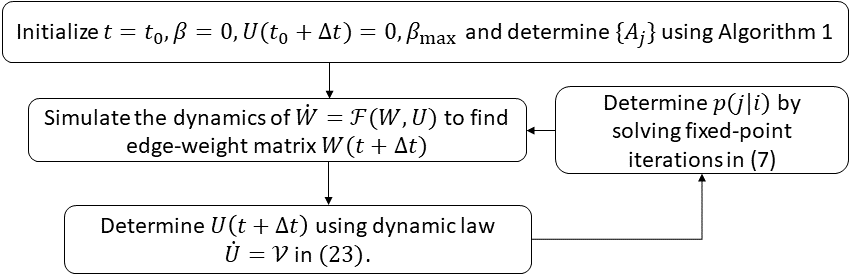}
    \caption{\small Algorithm to solve the multiway $k$-cut problem on Dynamic graphs.}\label{fig:flowchart}
\end{figure}
\begin{figure}
    \centering
    \includegraphics[height=2.2in,width=2.5in]{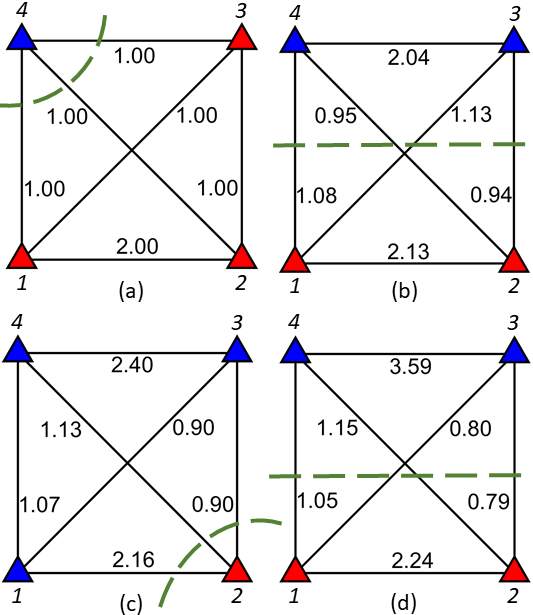}
    \caption{{\small Illustration of the proposed approach for time varying minimum multiway $2$-cut in an example graph with $4$ nodes, where $s_1=2$ and $s_2=4$. The four plots indicate different stages of the evolution of edge-weights. Note that the edge $(1,2)$ is uncontrollable and evolves under its natural dynamics. On the other hand, other edges can be controlled to result in minimum cut values at each instant without having to explicitly compute cut values from the beginning. It is easy to verify that the cut values at each instant.} }\label{fig:dyn4}
\end{figure}

{\em Dynamic Graphs: } Fig. \ref{fig:flowchart} summarizes the steps of our algorithm. Fig. \ref{fig:dyn4} illustrates the time-varying minimum multiway $2$-cut obtained on a weighted {\em undirected} dynamic graph using our proposed methodology. As stated in Section \ref{sec:DynamicCase} the objective here is to find the time varying cuts $A_1$ and $A_2$ where the nodes $s_1$ and $s_2$ are the nodes $2$ and $4$ respectively. The edge-weight dynamics is given  by $\dot{W}(t) = AW+BU$, which can be succinctly re-written as
\begin{align}\label{eq:dyn3_eqn}
{\small
\begin{pmatrix}
\dot{w}(1,2)\\
\dot{w}(1,3)\\
\dot{w}(1,4)\\
\dot{w}(2,3)\\
\dot{w}(2,4)\\
\dot{w}(3,4)
\end{pmatrix}=A'
\begin{pmatrix}
w(1,2)\\
w(1,3)\\
w(1,4)\\
w(2,3)\\
w(2,4)\\
w(3,4)
\end{pmatrix}+B'
\begin{pmatrix}
u(1,2)\\
u(1,3)\\
u(1,4)\\
u(2,3)\\
u(2,4)\\
u(3,4)
\end{pmatrix}, \text{where} }
\end{align}
\begin{align}\label{eq:dyn1_eq_par}
A'&=\frac{5}{6}\text{diag}\big(e^{0.5},~e^{1.8},~e^{1.5},~e^{-1.5},~e^{-1},~e^{3.5}\big)\nonumber, \text{ and }\\
B'&=\text{diag}\big(0,~1,~1,~1,~1,~1\big).
\end{align}
Note that in the above system dynamics (\ref{eq:dyn1_eq_par}), the edge $(1,2)$ is incapable of being influenced by any external control input. We simulate the dynamical system (\ref{eq:dyn3_eqn}) for a total of $0.075$ seconds. Fig. \ref{fig:dyn4}(a) illustrates the minimum multiway $2$-cut given by the Algorithm \ref{alg:alg1} at the initial time $t_0$ where $A_1 = \{1,2,3\}$, $A_2=\{4\}$. As time progresses the partitioning changes to $A_1=\{1,2\}$,$A_2=\{3,4\}$ in Fig. \ref{fig:dyn4}(b), $A_1=\{2\}$, $A_2=\{1,3,4\}$ in Fig. \ref{fig:dyn4}(c) and $A_1=\{1,2\}$, $A_2=\{3,4\}$ in Fig. \ref{fig:dyn4}(d). 

In the frame-by-frame approach we first discretize the entire time interval of $0.075$ seconds with $\Delta t = 0.01$ seconds and set $U(t)=0$ $\forall$ $t$, then Algorithm \ref{alg:alg1} is used to obtain the minimum multiway $2$-cut at each time instant. However this latter takes approximately $29$ times more computational time than required by our proposed method thereby making the frame-by-frame approach unscalable even for a very short duration of time.

\section{Conclusion}\label{sec:conclusion}
In this paper, we formulate a statistical physics based algorithm for the minimum multiway $k$-cut problem on static as well as dynamic digraphs and present a convergence proof for the algorithm. The algorithm described in this paper is computationally efficient and has ability to avoid poor local minima through controlled randomness. We believe that a combination of good theoretical properties and experimental success of the proposed MEP-based algorithm makes it a suitable technique of choice for a wide variety of combinatorial optimization problems on graphs, such as, vertex coloring and finding independent sets.

\section*{APPENDIX}\label{sec:appendix}
\subsection*{Proof of Theorem 1}
1): The first part of $L_1(W(t))$ in (\ref{eq:L1}) is positive since $\sum_{j}p(j|l)p(j|m)\leq (\sum_j p(j|l)^2)^{0.5}(\sum_j p(j|m))^{0.5}\leq (\sum_j p(j|l))^{0.5}(\sum_j p(j|m))^{0.5}=1$, and the second part $\sum_{ij}p(j|i)\log p(j|i)\geq -N\log k$. Also since $\|U\|_{\mathrm{F}}\geq 0$ we have that $F+\frac{1}{\beta}N\log k\geq 0$.

2(a): Since $\bar{\mathcal{V}}$ is Lipschitz at $U=0$, $\exists$ a $\delta$ neighborhood $B_{\delta}\triangleq\{U:\|U\|\leq \delta\}$ and $k>0$ such that $\dot{F}(\bar{\mathcal{V}}) = \frac{1}{2}\alpha + 2\mu e_N^T U\circ\bar{\mathcal{V}}e_N\leq 0$, where $\alpha = 2e_N^T\Phi\circ W\circ \dot{W}e_{N}$, and $\|\bar{\mathcal{V}}\|\leq k\|U\|$ $\forall$ $U\in B_{\delta}$.

\textit{Case $\alpha>0$: } $\small \Rightarrow |\alpha|<|4\mu e_N^TU\circ\bar{\mathcal{V}}e_N|\leq \bar{k}\|U\|_\mathrm{F}^2$. From $\mathcal{V}$ in (\ref{eq:cont_dyn}), we have
\begin{small}
\begin{align*}
    \|\mathcal{V}\|_\mathrm{F}&\leq \Big[C_0 + \frac{\bar{k}}{4\mu}+\frac{\sqrt{\bar{k}^2+16\mu^2}}{4\mu}\Big]\|U\|_\mathrm{F}= \hat{k}\|U\|_\mathrm{F}
\end{align*}
\end{small}
\noindent\textit{Case $\alpha <0$: } From (\ref{eq:cont_dyn})
\begin{small}
\begin{align*}
\|\mathcal{V}\|_\mathrm{F} \leq (C_0+\frac{\sqrt{\bar{k}^2+16\mu^2}}{4\mu})\|U\|_\mathrm{F}=\check{k}\|U\|_\mathrm{F}
\end{align*}
\end{small}
2(b): Substituting the dynamic controller $\mathcal{V}(t)$ in (\ref{eq:cont_dyn}) in the expression of $\dot{F}$ given by (\ref{eq:Vdot}) we obtain $\dot{F} = -2\mu C_0 \|U(t)\|_\mathrm{F}^2-\sqrt{\alpha^2+(4\mu \|U(t)\|_\mathrm{F}^2)^2}$ which is clearly non-positive.

2(c): $\dot{F}\rightarrow 0$ follows from Lasalle's Invariance Principle. From equation (\ref{eq:cont_dyn}) $\frac{d(U^TU)}{dt}\leq -C_0 U^TU$ which implies $U^TU \rightarrow 0$ exponentially using the Gronwall's inequality \cite{gronwall1919note}. This in turn implies that $U\rightarrow 0$. Now since $\dot{L}_1 = \dot{F} - 2e_N^T\mu U(t)\circ \mathcal{V}e_N$ we have that $\dot{L}_1\rightarrow 0$.

\bibliographystyle{IEEEtran}
\bibliography{IEEEabrv,mybibfile}

\end{document}